\documentclass[11pt,oneside]{amsart}

\usepackage{amsmath,ifthen, amsfonts, amssymb}
\usepackage{graphicx}

\newcommand{\showcomments}{no}

\newsavebox{\commentbox}
\newenvironment{com}%
{\ifthenelse{\equal{\showcomments}{yes}}%
{\footnotemark
    \begin{lrbox}{\commentbox}
    \begin{minipage}[t]{1.25in}\raggedright\sffamily\tiny
    \footnotemark[\arabic{footnote}]}
{\begin{lrbox}{\commentbox}}}%
{\ifthenelse{\equal{\showcomments}{yes}}%
{\end{minipage}\end{lrbox}\marginpar{\usebox{\commentbox}}}
{\end{lrbox}}}

\newtheorem{thm}{Theorem}[section]
\newtheorem{lem}[thm]{Lemma}
\newtheorem{cor}[thm]{Corollary}

\newtheorem{prop}[thm]{Proposition}
\newtheorem*{mainTRIPS}{Theorem~\ref{thm:TRips}}
\newtheorem*{applicationTOUT}{Theorem~\ref{thm:Tout}}
\newtheorem*{thmTnonHopf}{Theorem~\ref{thm:TnonHopf}}
\newtheorem*{thmTnoncoHopf}{Theorem~\ref{thm:TnoncoHopf}}

\theoremstyle{definition}
\newtheorem{defn}[thm]{Definition}

\newtheorem{rem}[thm]{Remark}

 \DeclareMathOperator{\Aut}{Aut}
\DeclareMathOperator{\Out}{Out} \DeclareMathOperator{\Inn}{Inn}

\DeclareMathOperator{\girth}{girth}
\DeclareMathOperator{\diam}{Diameter}

\newcommand{\scname}[1]{\text{\sf #1}}
\newcommand{\area}{\scname{Area}}

\newcommand{\field}[1]{\mathbb{#1}}
\newcommand{\integers}{\ensuremath{\field{Z}}}
\newcommand{\naturals}{\ensuremath{\field{N}}}

\newcommand{\boundary}   {{\ensuremath \partial}}

\renewcommand{\geq}{\geqslant}
\renewcommand{\leq}{\leqslant}
\renewcommand{\phi}{\varphi}
\renewcommand{\epsilon}{\varepsilon}

\newcommand{\norm}[1]{\left\|#1\right\|}
\newcommand{\abs}[1]{\left|#1\right|}

\newcommand{\Gr}{\mathit{Gr}'}

\begin{document}

\title[Kazhdan groups with infinite outer automorphism group]
{Kazhdan groups with infinite outer automorphism group}

\author[Y.~Ollivier]{Yann Ollivier}
      \address{CNRS, UMPA\\
      \'Ecole normale sup\'erieure de Lyon\\
      46, all\'ee d'Italie\\
      69364 Lyon cedex 7\\
                              France }
      \email{yann.ollivier@normalesup.org}

\author[D.~T.~Wise]{Daniel T. Wise}
      \address{Dept.\ of Math.\\
               McGill University \\
               Montr\'eal, Qu\'ebec, Canada H3A 2K6 }
      \email{wise@math.mcgill.ca}
\thanks{Research supported by NSERC grant}

\date{\today}

\begin{abstract}
For each countable group $Q$ we produce a short exact sequence
$1\rightarrow N \rightarrow G \rightarrow Q\rightarrow 1$ where
$G$ is f.g.\ and has a graphical $\frac16$ presentation and $N$ is
f.g.\ and satisfies property~$T$.

As a consequence we produce a group $N$ with property~$T$ such
that $\Out(N)$ is infinite.

Using the tools developed we are also able to produce examples
of nonHopfian and non-coHopfian groups with property~$T$.

One of our main tools is the use of random groups to achieve
certain properties.
  \end{abstract}

\maketitle

\section{Introduction}

The main result of this paper is a variant of Rips' construction
which allows us to get groups with infinite outer automorphism
group, combined with a tool of Gromov to get property $T$ at the same
time. Yet other variants provide non-Hopfian and non-coHopfian groups
with property $T$.

In \cite{Rips82}, Rips gave an elementary construction which given
a countable group~$Q$ produces a short exact sequence
$1\rightarrow N\rightarrow G\rightarrow Q\rightarrow 1$, where $G$
is a $C'(\frac16)$ group and $N$ is finitely generated. Rips used
his construction to produce $C'(\frac16)$ presentations with
various interesting properties, by lifting pathologies in $Q$ to
suitably reinterpreted pathologies in $G$.

Besides, Gromov~\cite{Gromov2003} was able to produce (random) groups
with property $T$ having so-called \emph{graphical~$1/6$ small
cancellation} (or $\Gr(\frac16)$ for short) presentation, which is a kind
of generalized $C'(\frac16)$ small cancellation property (see 
section~\ref{section:graphical} below). 

The mixture of these two tools yields the following
(section~\ref{section:TRips}):
\begin{thm}
\label{thm:TRips}
 For each countable group $Q$, there is a short
exact sequence $1\rightarrow N \rightarrow G \rightarrow
Q\rightarrow 1$ such that \begin{enumerate}
\item $G$ is torsion-free,
 \item $G$ has a graphical~$\frac16$ presentation, and
\item $N$ has property~$T$.
\item Moreover, $G$ is finitely generated if $Q$ is, and
finitely presented if $Q$ is.
\end{enumerate}
\end{thm}

The graphical~$1/6$ presentation keeps enough properties of ordinary
small cancellation as to mix nicely with Rips' construction. However,
we note that Theorem~\ref{thm:TRips} cannot be
obtained with $G$ an ordinary $C'(\frac16)$ group, since finitely
presented $C'(\frac16)$ groups act properly on a $C\!A\!T(0)$ cube
complex by~\cite{WiseSmallCanCube04}, and hence their infinite
subgroups cannot have Property~$T$~\cite{NibloReeves97,NibloRoller98}.

\begin{com}
There is a more subtle point here:
If $N$ has T, then there is $M\twoheadrightarrow N$ s.t. M
is f.p. with T, and since $N \subset G$, we see that $M$ factors through,
and so can be chosen to be a f.p. $C'(\frac16)$ group.
Thus the construction cannot work using ordinary $C'(\frac16)$ theory.
\end{com}

We apply Theorem~\ref{thm:TRips}
to obtain the following (section~\ref{section:Tout}):
\begin{thm}\label{thm:Tout}
There exists a group $N$ with property~$T$ such that $\Out(N)$ is
infinite.
\end{thm}

In fact we even prove that any countable group embeds in $\Out(N)$
for some Kazhdan group $N$.

 The motivation is that, as proven by
 Paulin \cite{Paulin91b},
 \begin{com}
 And Bestvina? Reference?
 \end{com}
  if $H$ is word-hyperbolic and
$|\Out(H)|=\infty$ then $H$ splits over an infinite cyclic group,
and hence $H$ cannot have property~$T$. The question of whether
every group with property~$T$ has a finite outer automorphism
group belongs to the list of open problems mentioned in de la Harpe and
Valette's classical book on Property $T$ (\cite{HarpeValette89}, p.~134),
was raised again by  Alain Valette in his mathscinet review of \cite{Paulin91b},
 and later appeared in a problem list from the 2002~meeting on
property~$T$ at Oberwolfach.

Finally, we use the tools we developed to obtain the following two
examples (sections~\ref{section:TnonHopf} and~\ref{section:TnoncoHopf}):
\begin{thm}\label{thm:TnonHopf}
There exists a Kazhdan group $G$ that is not Hopfian.
\end{thm}
\begin{thm}\label{thm:TnoncoHopf}
There exists a Kazhdan group $G$ that is not coHopfian.
\end{thm}

Various other attempts to
augment Rips's construction have focused on strengthening the
properties of $G$ when $Q$ is f.p.\ (e.g.:
$G$ is $\pi_1$ of a negatively curved complex \cite{Wise98};
$G$ is a residually finite $C'(\frac16)$ group \cite{WiseRFRips};
$G$ is a subgroup of a right-angled Artin group, so $G\subset
S\!L_n(\integers)$ \cite{HaglundWiseSpecial}).

One key ingredient of our constructions is the use of random methods,
introduced by Gromov~\cite{Gromov93} (see also~\cite{Ghys03}
and~\cite{Ollivier04} for a discussion of random groups), to provide
examples of groups with particular properties. Namely, we use a result
of~\cite{Gromov2003} providing a presentation of a group with property $T$
satisyfing the graphical small cancellation property. We include in
section~\ref{section:Gromov graph} a standalone proof of the results we need
from~\cite{Gromov2003}.

\section{$\Gr(\frac16)$ graphs}
\label{section:graphical}

\subsection{Review of graphical $\alpha$-condition $\Gr(\alpha)$}
Throughout all this article, $B$ is a bouquet of $m\geq 2$ circles whose
edges are directed and labelled, so that $m$ will be the number of
generators of the group presentatations we consider.

Let $\Gamma\looparrowright B$ be an immersed graph,
and note that $\Gamma$ has an induced labelling. That $\Gamma$
immerses in $B$ simply denotes the fact that the words carried by
paths immersed in $\Gamma$ are reduced.

By definition, the group $G$ presented by $\langle B | \Gamma
\rangle$, has generators the letters appearing on $B$, and
relations consisting of all cycles appearing in $F$.

A {\em piece} $P$ in $\Gamma$ is an immersed path
$P\looparrowright B$ which lifts to $\Gamma$ in more than one way.

\begin{defn}
We say $\Gamma\looparrowright B$ satisfies the {\em graphical
$\alpha$ condition} $\Gr(\alpha)$ if for each piece~$P$,
 and each cycle $C\rightarrow \Gamma$ such that $P\rightarrow \Gamma$ factors
 through $P\rightarrow C\rightarrow \Gamma$,
we have $|P| < \alpha |C|$.
\end{defn}

\begin{com}
PICTURE here featuring a piece in bold appearing on two cycles,
and some other, maybe shorter cycles.
\end{com}

The graphical $\alpha$ condition generalizes the usual
$C'(\alpha)$: let $F$ consist of the disjoint union of a set of
cycles corresponding to the relators in a presentation. The
graphical $\alpha$ condition is a case of a complicated but more
general condition given by Gromov~\cite{Gromov2003}.

The condition $\Gr(\frac16)$ implies that the group $G$ is torsion-free,
word-hyperbolic, of dimension $2$, just as the
$C'(1/6)$ condition \cite{Ollivier03}. The group is non-elementary except
in some explicit degenerate cases (a hyperbolic group is called
elementary if it is finite or virtually $\integers$).

There is also a slightly stronger version of this condition, in which we
demand that the size of the pieces be bounded not by $\alpha$ times the
size of any cycle containing the piece, but by $\alpha$ times the girth
of $\Gamma$ (recall the \emph{girth} of a graph is the smallest length of
a non-trivial closed path in it). We will sometimes directly prove this
stronger version below, since it allows lighter notations.

A \emph{disc van Kampen diagram} w.r.t.\ a graphical presentation is a van
Kampen diagram every $2$-cell of which is labelled by a closed path immersed in
$\Gamma$. It is \emph{reduced} if, first, it is reduced in the ordinary
sense and if moreover, for any two adjacent $2$-cells, the boundary word
of their union does not embed as a closed path in $\Gamma$ (otherwise,
these two $2$-cells can be replaced by a single one). It is proven
in~\cite{Ollivier03} that if $\Gr(\alpha)$ holds, such a reduced van
Kampen diagram satisfies the ordinary $C'(\alpha)$ condition.

\subsection{Producing more $\Gr(\frac16)$ graphs}
One useful feature of a presentation satisfying
the ordinary $C'(\frac16$) theory
is that, provided that the relations are not ``too dense'' in
a certain sense, more relations can be added to the presentation
without violating the $C'(\frac16)$ condition.

In this subsection, we describe conditions on a $\Gr(\alpha)$
presentation such that additional relations can be added.

\begin{prop}\label{prop:adding cycles}
Let $\Gamma\looparrowright B$ satisfy the $\Gr(\alpha)$ condition
and suppose there is an immersed path $W\rightarrow B$ such that
$1\leq|W|<\frac\alpha{2}\girth(\Gamma)-1$, and $W$ does not lift
to $\Gamma$.

Then there is a set of closed immersed paths $C_i\looparrowright
B: i\in \naturals$ such that the disjoint union
$\Gamma'=\Gamma\sqcup_{i\in\naturals} C_i\looparrowright B$
satisfies the $\Gr(\alpha)$ condition.
\end{prop}

\begin{proof}

We first form an immersed labelled graph $A\looparrowright B$ as
follows: Let $D$ be the radius~$2$ ball at the basepoint of the
universal cover $\tilde B$, and attach two copies $W_x$ and $W_y$
 of the arc $W$ along four distinct leaves of $D$
 as in Figure~\ref{fig:Agraph}.
 (This can always be done avoiding the inverses of the initial and final letter of $W$, so that $D$ immerses in $B$).
Finally, we remove the finite trees that remain.

\begin{figure}\centering
\includegraphics[width=.5\textwidth]{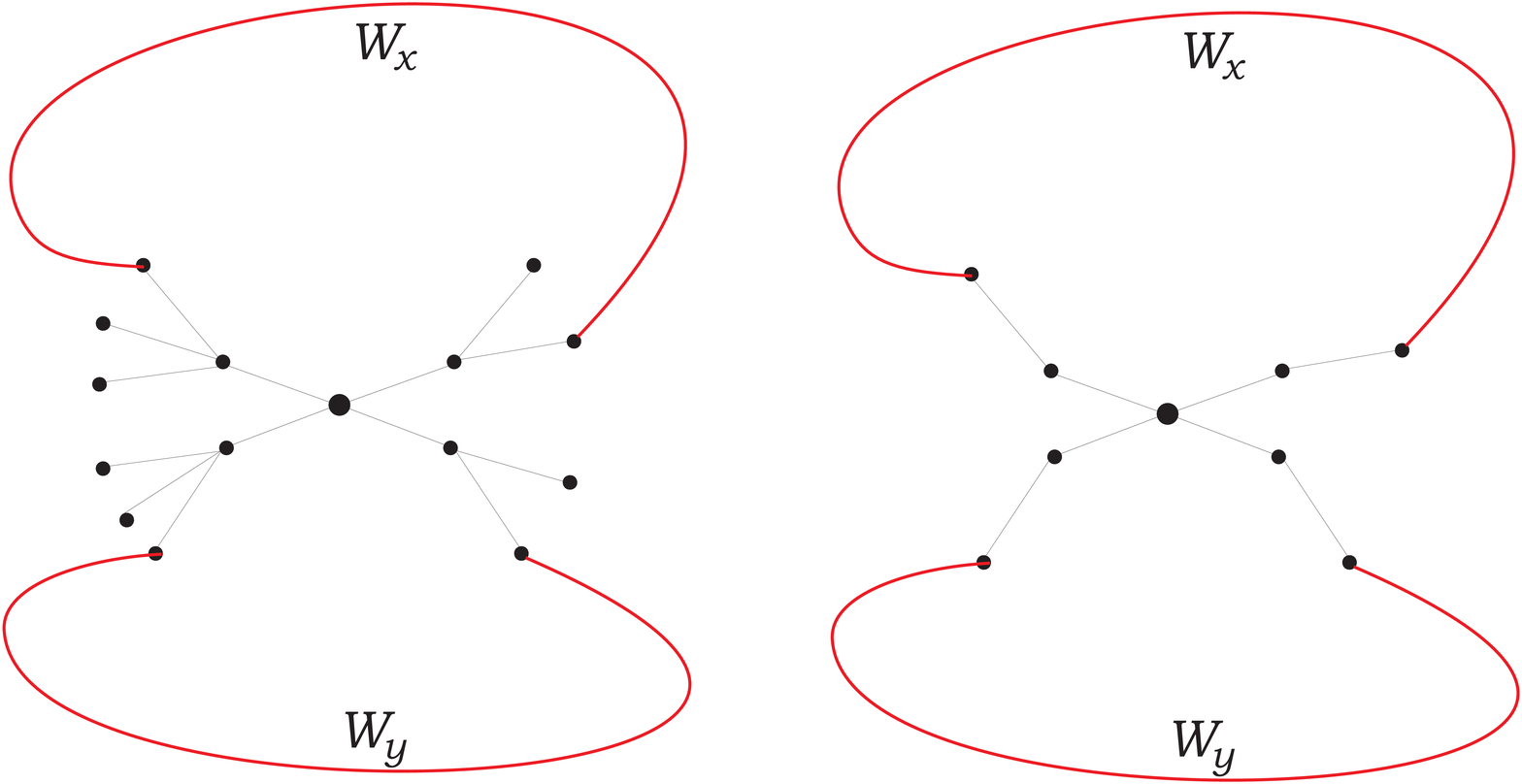}
\caption{ \label{fig:Agraph}}
\end{figure}
\begin{com}
There are some missing edges in the left figure: the degree should be 4 at
each point, not 3.
\end{com}

Observe that any path $P\looparrowright B$ that lifts to both $A$ and
$\Gamma$ satisfies $|P|< \alpha\girth(\Gamma)$. Indeed, if $P$
lifts to $\Gamma$, then $P$ cannot contain $W_x$ or $W_y$ as a
subpath, and hence $P=U_1U_2U_3$ where $U_1$ and $U_3$ are proper
initial or terminal subpaths of a $W$-arc, and $U_2$ is a path in
$D$, so $|P|\leq |U_1|+|U_2|+|U_3| \leq (|W|-1) + 4 + (|W|-1) =
2|W|+2 =2(|W|+1)<
2\frac\alpha{2}\girth(\Gamma)=\alpha\girth(\Gamma)$.

Now let $x$ and $y$ be arbitrary labels.
To any reduced word $w$ in the letters $x^{\pm 1}$ and $y^{\pm 1}$ we can
associate an immersed closed path $\phi(w)$ in $A$ by sending $x$
to the based path in $A$ containing $W_x$, and similarly for $y$.

Now for each $i\in\naturals$, let
 $c_i$ denote the word
 $xy^{1000i+1}xy^{1000i+2}\cdots xy^{1000i+999}$.
 It is easily verified that for large enough values of $1000$, the set of words $\langle x,y\mid c_i: i\in\naturals\rangle$
satisfies the $C'(\frac\alpha{2})$ condition.

 Let $C_i\looparrowright B$
denote the corresponding closed immersed cycle $\phi(c_i)$. Pieces
in $\bigsqcup C_i$ are easily bounded in terms of pieces in $\langle
x,y\mid c_i \, (i\in \naturals)\rangle$, so that $\bigsqcup C_i$
satisfies the $\Gr(\alpha)$ (actually $C'(\alpha)$) condition.

Finally $\Gamma'= \Gamma \sqcup_{i\in \naturals} C_i$ satisfies
the $\Gr(\alpha)$ condition since pieces that lift twice to
$\Gamma$ are bounded by assumption, and we have just bounded
pieces that lift to $\Gamma$ and to some $C_i$, and pieces that
lift to some $C_i$ and some $C_j$.
\end{proof}

\begin{rem}
\label{rem:nontriv}
The missing word condition in $\Gamma$ ensures that the group presented
by $\Gamma$ is non-elementary. Indeed, the group presented by
$\Gamma\sqcup_i C_i$ has infinite Euler characteristic (it is of
dimension $2$) and is thus
non-elementary, so a fortiori the group presented by $\Gamma$ is.
\end{rem}

\section{The T Rips construction}
\label{section:TRips}

Let us now turn to the proof of the main theorem of this article.
We use an intermediate construction due to Gromov.

\begin{prop}
\label{prop:T graph}
There exists a finite graph $\Gamma$ that immerses in a bouquet
$B$ of two circles such that:
\begin{enumerate}
\item The group presented by $\langle B \mid \Gamma\rangle$ has
property~$T$.

\item $\Gamma\looparrowright B$ satisfies the $\Gr(\frac1{12})$ condition.

\item There is a path $W\looparrowright B$ with
$1\leq |W|<\frac{1}{24}\girth(\Gamma)-1$ and $W$ does not lift to
$\Gamma$.

\item $\Gamma$ has arbitrarily large girth.
\end{enumerate}
\end{prop}

A proof of this is included below (section~\ref{section:Gromov
graph}).

\begin{mainTRIPS}
For each countable group $Q$, there is a short exact sequence
$1\rightarrow N \rightarrow G \rightarrow Q\rightarrow 1$
such that \begin{enumerate}
\item $G$ is torsion-free,
 \item $G$ has a graphical~$\frac16$ presentation, and
\item $N$ has property~$T$ and is non-trivial.
\item Moreover, $G$ is finitely generated if $Q$ is, and finitely
presented if $Q$ is.
\end{enumerate}
\end{mainTRIPS}

\begin{proof}
Let $Q$ be given by the following presentation:
$$\langle q_i: i\in I \mid R_j: j\in J \rangle$$
Let $\Gamma\looparrowright B$ be a graph provided by
Proposition~\ref{prop:T graph}, where the
edges of $B$ are labelled by $x$ and $y$. Let
$\Gamma'=\Gamma\sqcup_n C_n$ be as in Proposition~\ref{prop:adding
cycles} with $\alpha=1/12$.

The presentation for $G$ will be the following:
\begin{multline}\label{pres:Trips}
\langle\ 
x,y,\, q_i \ (i \in I) \ \mid \ 
\Gamma, \\
 x^{q_i}=X_{i+},\ 
x^{q_i^{-1}}=X_{i-},\ 
y^{q_i}=Y_{i+},\ 
y^{q_i^{-1}}=Y_{i-} \ (i\in I),\\
R_j=W_j \ (j\in J)
\ \rangle
\end{multline}
where superscripts denote conjugation, and where the $X_{i+}$, $X_{i-}$,
$Y_{i+}$, $Y_{i-}$, and $W_j$
are equal to paths corresponding to distinct $C_n$ cycles of $\Gamma'$,
$|W_j|>12|R_j|$ for each $j\in J$,
and $|X_{i\pm}|>36, |Y_{i\pm}|>36$ for each $i\in I$.

The $\frac16$ condition follows easily. Let us check, for example,
that there is no $\frac16$-piece between $\Gamma$ and the
relation $x^{q_i}=X_{i+}$. Since the $q_i$'s do not appear as labels on $\Gamma$,
any such $\frac16$-piece would be either $x$ or a subword of $X_{i+}$.
The former is ruled out since $\girth(\Gamma)>6$. The latter would provide
a piece between $\Gamma$ and $X_{i+}$ (which is one of the $C_n$'s);
such a piece is by assumption of length at most $\frac1{12}|X_{i+}|$
which in turn is less than $\frac16|x^{q_i}=X_{i+}|$ as needed.
The other cases are treated similarly.

Now $N$ is the subgroup of $G$ generated by $x$ and $y$.
It is normal by construction of the presentation of $G$.
Note that $N$ has property~$T$ since it is a quotient of $\langle
x,y\mid \Gamma\rangle$ which has property~$T$ by choice of
$\Gamma$.

Finally, $N$ is non-trivial: indeed, we can pick some cycle $C_n$ which is a
word in $x,y$ and which will be in small cancellation with the rest of
the presentation. This provides a word in $x$ and $y$ whichh is not
trivial in the group.
\end{proof}

\section{Kazhdan groups with infinite outer automorphism group}
\label{section:Tout}

\begin{applicationTOUT}
Any countable infinite group $Q$ embeds
in $\Out(N)$ for some group $N$ with property~$T$.

In particular, there exists a group $N$ with property~$T$ such that $\Out(N)$ is
infinite.
\end{applicationTOUT}

\begin{proof}
For $1\rightarrow N \rightarrow G\rightarrow Q\rightarrow 1$, the
group $G$ acts by inner automorphisms on itself, so we have a
homomorphism $G\rightarrow \Aut(N)$, and $N$ obviously maps to
$\Inn(N)$ so there is an induced homomorphism $Q=G/N \rightarrow
\Out(N)$. Elements in the kernel of $Q\rightarrow \Out(N)$ are
represented by elements $g\in G$ such that $m^g=m^n$ for some
$n\in N$ and all $m\in N$. Thus $gn^{-1}$ centralizes $N$.

First suppose that $Q$ is finitely presented, so that $G$ is as well.

In this case $N$ is a non-elementary subgroup of the torsion-free
word-hyperbolic group $G$, and hence $N$ has a trivial centralizer.
Indeed, $N$ must contain a rank~$2$ free subgroup $\langle n_1,n_2\rangle$
(see~\cite{GhysBook90}, p.~157).
If a nontrivial element $c$ centralizes $N$ then $\langle c,n_1\rangle$
and $\langle c,n_2\rangle$ are both abelian, and hence infinite cyclic
since $G$ cannot contain a copy of $\integers^2$.
Thus $n_1^{m_1}=c^{p_1}$ and $n_2^{m_2}=c^{p_2}$
for some $p_i, m_i \neq 0$.
But then $n_1^{m_1}$ commutes with $n_1^{m_2}$ which is impossible.

Since the centralizer of $N$ is trivial,
we have $gn^{-1}=1$, so $g\in N$, and hence $Q\rightarrow
\Out(N)$ is injective.

The case when $Q$ is not finitely presented reduces back to the previous
one:
Indeed, suppose that some element $g$ of $G$ lies in the centralizer of $N$.
This is equivalent to stating that $g$ commutes with $x$ and $y$. But $g$ can
be written as a product of finitely many generators, and similarly the
relations $[g,x]=1$ and $[g,y]=1$ are consequences of only finitely many
relators, so that $g$ still lies in the centralizer of $N$ in a finite
subpresentation of the presentation of $G$.
\end{proof}

\begin{rem}
By adding some additional relations to $N$, the above argument was
used in \cite{BumaginWiseOut} to show that every countable group
$Q$ appears as $\Out(N)$ for some f.g.\ $N$, and that every
f.p.\ $Q$ appears as $\Out(N)$ where $N$ is f.g.\ and residually finite
(but property $T$ did not appear there).

It appears likely that a more careful analysis along those lines,
would show that every countable group arises as $\Out(N)$ where
$N$ has property~$T$.
\end{rem}

\section{A Kazhdan group that is not coHopfian}
\label{section:TnoncoHopf}

\begin{thmTnoncoHopf}
There exists a Kazhdan group that is not coHopfian.
\end{thmTnoncoHopf}

\begin{proof}
Consider the group
$$G=\langle\, a,b,t\mid \Gamma, \, a^t=\phi(a), \, b^t=\phi(b) \, \rangle$$
where $\phi(a)$ and $\phi(b)$ are chosen so that $\Gamma\sqcup
\phi(a)\sqcup \phi(b)$ satisfies $\Gr(\frac16)$ and $\abs{\phi(a)}>3$,
$\abs{\phi(b)}>3$. (This is in fact a
subpresentation
of the presentation~(\ref{pres:Trips}) used in the proof of
Theorem~\ref{thm:TRips}.)


Clearly, the subgroup $K=\langle a,b\rangle$ is a Kazhdan group since
it is a quotient of $\langle a,b\mid \Gamma\rangle$.

The map $K\rightarrow K$ induced by $\phi$ is clearly well-defined and
injective since it arises from conjugation in the larger group $G$.

We will now show that $\phi$ is not surjective by verifying that
$a\not\in\langle \phi(a),\phi(b)\rangle$.

We argue by contradiction: Suppose that $a$ is equal in $G$ to a word
$W(\phi(a),\phi(b))$ in $\phi(a)$ and $\phi(b)$;
we can choose $W$ such that the disc diagram expressing this equality
in the presentation for $G$ has minimal area among all such choices.
Note that since $D$ is reduced and $G$ is $\Gr(\frac16)$,
$D$ is a diagram satisfying the ordinary $C'(\frac16)$ condition.

By Greendlinger's Lemma, (after ignoring trees possibly
 attached to $\boundary D$) either $D$ is a single $2$-cell,
or $D$ has at least two $2$-cells whose outer paths are the majority
of their boundaries.

The first possibility is excluded by consideration of the presentation for $G$.
In the second case, one such $2$-cell $R$ has outerpath $Q$ not containing
the special $a$-edge in $\boundary D$,
as illustrated on the left in Figure~\ref{fig:tannulusremoval}.

\begin{figure}\centering
\includegraphics[width=\textwidth]{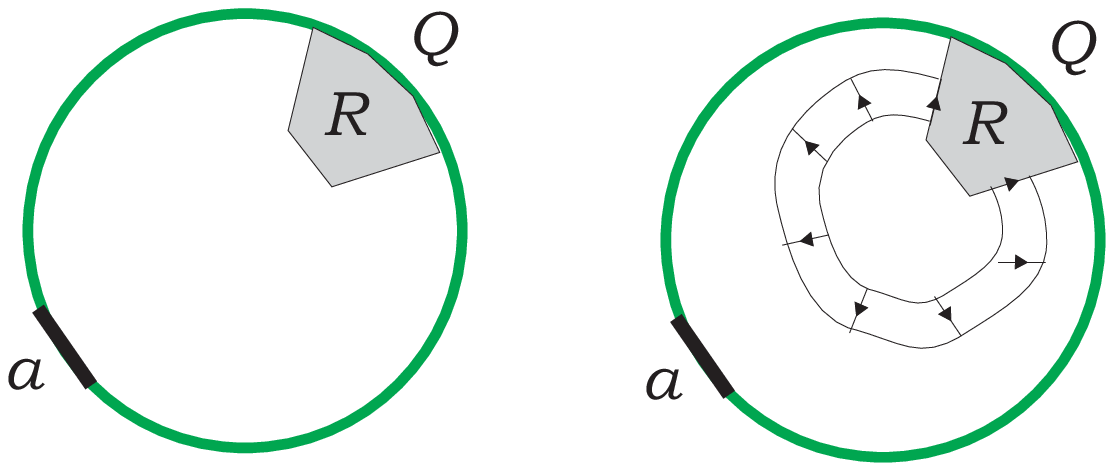}
\caption{ \label{fig:tannulusremoval}}
\end{figure}
\begin{com}
In the figure, should add the t's and shold make more obvious what is D'
\end{com}

The boundary word of $2$-cell $R$ cannot be a words immersing in
$\Gamma$.
Indeed, since it has more than half its length on the boundary of $D$ and
this boundary bears a word in $\phi(a)$ and $\phi(b)$, this would
contradict the small cancellation property of $\Gamma\sqcup
\phi(a)\sqcup\phi(b)$.
So $R$ is a $2$-cell expressing the equality $a^t=\phi(a)$ or
$b^t=\phi(b)$. Moreover, since $t$ does not appear on the boundary of
$D$, the side of $R$ on the boundary is the $\phi$-side.

Since $t\not\in \boundary D$, we can find a $t$-annulus
containing $R$ as illustrated in the center of in
Figure~\ref{fig:tannulusremoval}.

We now produce a new diagram $D'$ with $\area(D')<\area(D)$.
We do this by travelling around the $t$-annulus as on the right
in Figure~\ref{fig:tannulusremoval}.

Observe that the small cancellation property implies that an edge in the
$\phi(a)\subset \boundary R$ or $\phi(b) \subset \boundary R$ lines up
with an edge in some $\phi(a)$ or $\phi(b)$ in $\boundary D$, and at
exactly the same position. So if the $\phi(a)$ or $\phi(b)$ of $\boundary
R$ is not wholly contained in $\boundary D$, after removing $R$ the words
on the paths from $\boundary D$ to the $t$-edges of $R$ will cancel with
corresponding subwords of $\phi(a)$ and $\phi(b)$ lying in the remaining
part of $\boundary D$.

This implies that, after removing the annulus,
the boundary of $D'$ is labelled (maybe after folding) by a word of the
form $a=W'(\phi(a),\phi(b))$. But this is a contradiction since $D$ was
assumed to be minimal.

(Note that $D'$ might touch the special $a$-edge, and $D'$ might
have some extra singular edges.)
\end{proof}

\section{A Kazhdan group that is not Hopfian}
\label{section:TnonHopf}

\begin{defn}
Let $B$ be a bouquet of circles, and let $\phi:B\rightarrow B$.
Let $A\rightarrow B$ be a map of graphs,
then we let $\phi(A)\rightarrow B$ be the new map of graphs
where $\phi(A)$ is obtained from $A$ by substituting an arc $\phi(e)$
for each edge $e$ of $A$. That is, we replace the label on each edge of
$A$ by its image under $\phi$.
\end{defn}

\begin{lem}\label{lem:subsublemma}
Let $\Delta$ be a labelled graph satisfying $\Gr(\alpha)$ and
$\alpha\girth(\Delta)\geq1$. Suppose there is a path $P\looparrowright
\Delta$ such that the edges in $P$ all bear the same label $a$ and such
that $P$ factors through a closed path $P\looparrowright C\looparrowright
\Delta$. Then $\abs{P}< 2\alpha\abs{C}$.
\end{lem}

Note that the assumption $\girth(\Delta)\geq1/\alpha$ is not very strong:
if $\alpha\girth(\Delta)\leq1$ then a single letter can constitute a
piece, which can result in various oddities.
This lemma is false for trivial reasons if we remove this girth
assumption: when $\girth(\Delta)=1$ there are arbitrarily long
homogeneous paths though $\Gr(0)$ may be satisfied.

\begin{proof}
First, let us treat the trivial case when there is a length-$1$ loop
bearing label $a$: this implies $\girth(\Delta)=1$ so $\alpha=1$ and the
equality to show is trivial. The case $\abs{P}=1$ is trivial as well.

Second, suppose that there is no length-$1$ loop. Let $P$ be a path
labelled by $a^s$ with $s\geq 2$. Then the two paths labelled by
$a^{s-1}$ obtained by removing the first and last edge of $P$
respectively constitute a piece, and so we have $s-1<\alpha \abs{C}$ so
that $\abs{P}=s<\alpha\abs{C}+1\leq \alpha(\abs{C}+\girth(\Delta))\leq
2\alpha\abs{C}$.
\end{proof}

\begin{lem}\label{lem:sublemma}
Let $\Delta$ be a labelled graph satisfying $\Gr(\alpha)$ with
$\alpha\girth(\Delta)\geq 1$, and let $\phi:B\rightarrow B$
be induced by $a\mapsto a^n$ and $b\mapsto b^n$ for some $n\geq 1$.
Then, for any $k\in \naturals$, $\phi^k(\Delta)$ satisfies $\Gr(2\alpha)$.
\end{lem}

(Once more the girth assumption discards some degenerate cases when a
single edge can make a piece.)

\begin{proof}
The reader should think of $\phi^k(\Delta)$ as the $n^k$-subdivision of $\Delta$
where each $a$-edge is replaced by an arc of $n^k$ $a$-edges and likewise
for $b$-edges.

We begin by considering a homogeneous piece $P=a^r$ (or $P=b^r$ which is
similar) occurring in some cycle $C\looparrowright\phi^k(\Delta)$.  Then
$P$ is a subpath of a path $\phi^k(P')$ where $P'=a^{r'}$ is a path in
$\Delta$ and $P'$ occurs in a cycle $C'$ corresponding to $C$.

By the previous lemma, $|P'|=r'< 2\alpha|C'|$ and so
$|P|=r<2\alpha n^k|C'| = 2\alpha|C|$.

We now consider the general case where $P$ contains both $a$ and $b$ letters.
We may assume that $P$ is a maximal piece, in which case
$P=W(a^{n^k}, b^{n^k})$ where $P'=W(a,b)$ is itself a corresponding piece in $\Delta$.
Everything scales by $n^k$ i.e.
$|P|=n^k|P'|<n^k\alpha|C'|=\alpha|C|.$
\end{proof}

\begin{lem}\label{lem:cup phi n delta}
Let $\Delta$ satisfy $\Gr(\alpha)$ and suppose that
$\girth(\Delta)>1/\alpha$.
Let $n$ satisfy $n>s$ where $s$ is the maximal length of a path $a^s$ or
$b^s$ lifting to $\Delta$.
Let $\phi:B\rightarrow B$ be induced by $a\mapsto a^n$ and $b\mapsto b^n$.

Then $\bigsqcup_{k\geq 0}\, \phi^k(\Delta)$ satisfies $\Gr(8\alpha)$.
\end{lem}

Note that $s=\infty$ implies either $\girth(\Delta)=1$ (which is excluded
by assumption) or $\alpha=1$ (by removing the first and last letter of an
arbitrarily long $a^s$-path) in which case the affirmation is void. So
we can suppose $s<\infty$.

\begin{proof}
First, by the previous lemma, each $\phi^k(\Delta)$ itself satisfies $\Gr(2\alpha)$.

We now consider a piece $P$ between $\Delta$ and $\phi^k(\Delta)$.
Either $P\looparrowright \phi^k(\Delta)$ is contained in two subdivided
edges of $\phi^k(\Delta)$ so $|P|<2n^k$;
or $P$ contains an entire subdivided edge and hence an
$a^{n^k}$ (or $b^{n^k}$) subpath.

In the latter case when $P$ contains an $a^{n^k}$ or $b^{n^k}$ subpath,
since $P\looparrowright\Delta$ is a path in $\Delta$ then
$n^k$ is at most the maximal length of an $a$-path or $b$-path in $\Delta$.
But by hypothesis on $n$, this maximal length is bounded by~$n$,
and so $n^k<n$ which is impossible for $k\geq 1$.

In the former case, $P$ is the product of at most two homogeneous paths
(i.e.\ $a$-paths or $b$-paths) one of which has length~$\geq \frac12|P|$.
Thus by Lemma~\ref{lem:subsublemma}, $\frac12|P|<2\alpha|C|$ for any
cycle $C$ in $\Delta$ containing $P$. So $|P|<4\alpha|C|$ and so $P$
cannot be a $4\alpha$-piece in $\Delta$. Besides, suppose that $P$ is
included in a cycle $C$ immersed in $\phi^k(\Delta)$. Since
$\abs{P}<2n^k$ and $\abs{C}\geq
\girth{\phi^k(\Delta)}=n^k\girth(\Delta)\geq n^k/\alpha$ by assumption,
$P$ cannot consitute a $2\alpha$-piece in $\phi^k(\Delta)$ either.  (Note
that the constant $4$ is this reasoning is optimal: consider for $\Delta$
a circle of length $100$ containing $aabb$ at one place, and some garbage
for the rest; take $\alpha=(1+\epsilon)/100$ so that the two $a$'s do not
form a piece. Then $\phi(\Delta)$ contains some $aabb$ as well, so that
this word constitutes a $4/100$-piece in $\Delta\sqcup \phi(\Delta)$.)

Finally, we consider pieces between
$\phi^k(\Delta)$ and $\phi^{k'}(\Delta)$ where we can suppose $k'>k$.
We have just proved that $\Delta\sqcup \phi^{k'-k}(\Delta)$ satisfies $\Gr(4\alpha)$.
We now apply Lemma~\ref{lem:sublemma} to see that
$\phi^k(\Delta)\sqcup
\phi^{k'}(\Delta)=\phi^k\left(\Delta\sqcup\phi^{k'-k}(\Delta)\right)$
satisfies $\Gr(8\alpha)$.
\end{proof}

\begin{rem}
A generalization of Lemma~\ref{lem:cup phi n delta} should
hold with $\phi(a)$ and $\phi(b)$ appropriate small cancellation words
instead of $a^n$ and $b^n$.
\end{rem}

\begin{thmTnonHopf}
There exists a Kazhdan group that is not Hopfian.
\end{thmTnonHopf}
\begin{proof}
Let $G$ have the following presentation:
$$\langle\, a,b\mid \phi^i(\Gamma), \, \phi^i(a\phi(C_1)),\,
\phi^i(b\phi(C_2)), \, \phi^i(\phi(C_3)) \  (i\geq 0) \, \rangle$$
where
\begin{enumerate}
\item $\Gamma\sqcup C_1\sqcup C_2\sqcup C_3$ satisfies the $\Gr(\alpha)$
condition with $\alpha=1/2000$ ($C_1,C_2$ and $C_3$ arise
from Proposition~\ref{prop:adding cycles});
\item $\phi$ is defined by $\phi(a)=a^n$ and $\phi(b)=b^n$, for some $n$
greater than the maximal length of an $a$-word or $b$-word
in $\Gamma\sqcup C_1 \sqcup C_2\sqcup C_3$;
\item $\girth(\Gamma\sqcup C_1\sqcup C_2\sqcup C_3)\geq 2000$.
\end{enumerate}

Let $\Delta_0=\bigsqcup_{k\geq 0} \,\phi^k\left(\Gamma\sqcup C_1 \sqcup
C_2\sqcup C_3\right)$. By Lemma~\ref{lem:cup phi n delta}, this labelled
graph satisfies $\Gr(8\alpha)$.  As a subgraph of $\Delta_0$, the graph
$\Delta=\Gamma \sqcup\phi(C_1)\sqcup\phi(C_2)\sqcup C_3$ satisfies
$\Gr(8\alpha)$ as well.

We now prove that $\Delta'=\Gamma\sqcup a\phi(C_1) \sqcup b\phi(C_2)
\sqcup C_3$ is $\Gr(26\alpha)$.  Let $P$ be a piece involving the new
$a$-edge or the new $b$-edge.  Observe that $P=P_1aP_2$ (or $P_1bP_2$).
Note that a new $b$ (or new $a$) may lie in at most one of of $P_1$ or
$P_2$.  Thus $P$ is the concatenation of at most $3$~pieces in $\Delta$
together with the new~$a$ and possibly the new $b$.  Consequently for any
cycle $C$ containing $P$ we have $|P|<24\alpha|C|+2\leq26\alpha|C'|$
where we have used the hypothesis that $\alpha\girth\geq1$.

We now apply Lemma~\ref{lem:cup phi n delta} to see that
$\Omega=\bigsqcup_{k\geq 0} \phi^k(\Delta')$ satisfies $\Gr(208\alpha)$,
and so does the presentation for $G$ which is a subset of $\Omega$.

Now $\phi$ obviously sends relations to relations and thus induces a
well-defined map in $G$. This map is surjective since $a\phi(C_1)=_G 1$
and $b\phi(C_2)=_G 1$.

Finally $\phi$ is not injective since 
$\phi(C_3)=_G 1$ but $C_3\neq_G 1$.
Indeed, $C_3$ is in small cancellation relative to the relators of $G$
since both are included in $\Omega$.
\end{proof}

\section{A $T$ $\Gr(\frac16)$ graph with a missing word}
\label{section:Gromov graph}

A main point in this paper is the following, introduced by Gromov in~\cite{Gromov2003}:
\begin{prop}
\label{prop:sixth T W} For each $\alpha>0$ and $\alpha'>0$ there
exists a finite graph $\Gamma$ that immerses in a bouquet $B$ of
two circles such that:
\begin{enumerate}
\item The group presented by $\langle\, B \mid \Gamma\,\rangle$ has
property~$T$.

\item $\Gamma\looparrowright B$ satisfies the $\Gr(\alpha)$ condition.

\item There is a path $W\looparrowright B$ with
$1\leq|W|\leq\alpha'\girth(\Gamma)$ and $W$ does not lift to $\Gamma$.
\end{enumerate}

Moreover, the girth of $\Gamma$ can be taken arbitrarily large.
\end{prop}

This trivially implies Proposition~\ref{prop:T graph}. It also results
from Remark~\ref{rem:nontriv} that the obtained group is non-trivial.

The goal of the introduction of such graphs in~\cite{Gromov2003} was to
construct a group whose Cayley graph contains a family of expanders, in
relation with the Baum-Connes conjecture (see also~\cite{Ghys03}
and~\cite{Ollivier03expanders}). There, the construction is done starting
not only with a free group but with an arbitrary hyperbolic group
(compare~\cite{Ollivier04}), so that it can be iterated in order to embed
a whole family of graphs.

Here we use this construction for purposes closer to combinatorial group
theory. We do not need the full strength of the iterated construction;
this section is devoted to the proof of the statements we need.

We will use the following fact, the credit of which can be shared between
Lubotzky, Margulis, Phillips, Sarnack, Selberg. We refer
to~\cite{Lubotzky94} (Theorem~7.4.4 referring to Theorem~7.3.12), or to
to~\cite{DavidoffSarnackValette03}.

\begin{prop}\label{prop:expander fact}
For lost of $v\in \naturals$,
there is a family of graphs $\Gamma_i:i\in\naturals$ such that the
following hold:
\begin{enumerate}
\item Each $\Gamma_i$ is regular of valence $v$.

\item $\inf_i \lambda_1(\Gamma_i)>0$ where $\lambda_1$ denotes
the smallest non-zero eigenvalue of the discrete Laplacian
$\Delta$.

\item $\girth(\Gamma_i)\longrightarrow\infty$.

\item $\exists C$ such that $\diam(\Gamma_i)\leq
C\girth(\Gamma_i)$ for all $i$.
\end{enumerate}
\end{prop}

``Lots of $v$'' means e.g.\ that this works at least for $v=p+1$ with
$p\geq 3$ prime (\cite{Lubotzky94}, paragraph 1.2 refers to other
constructions). This is irrelevant for our purpose. 

We are going to use random labellings of subdivisions of the graphs $\Gamma_i$.
Subdividing amounts to labelling each edge with a long word rather
than just one letter, so that the small cancellation condition is
more easily satisfied.

That the diameter of the graph is bounded by a
constant times the girth reflects the fact that there are ``not too
many'' relations added (compare the density model of random groups
in~\cite{Gromov93} or~\cite{Ollivier04}): this amounts to taking
an arbitrarily small density.

To prove Proposition~\ref{prop:sixth T W} we need two more
propositions.

\begin{prop}
\label{prop:random T} Given $v\in\naturals$, $\lambda_0>0$ and an integer $j\geq 1$
there exists an explicit $g_0$ such that if $\Gamma$ is a $v$-regular graph
with $\girth(\Gamma)\geq g$, $\lambda_1(\Gamma)\geq \lambda_0$ and
$\Gamma$ is trivalent, then the random group defined through a
random labelling of the $j$-subdivision $\Gamma^j$ of $\Gamma$
will have property $T$, with probability tending to $1$ as the
size of $\Gamma$ tends to infinity.
\end{prop}

This is proven in \cite{Silberman03} (Corollary~3.19 where $d$ is our
$v$, $k$ is our number of generators $m$, and $\abs{V}$ the size of the
graph; in this reference, $\lambda(\Gamma)$ denotes the largest
eigenvalue not equal to $1$ of the averaging operator $1-\Delta$, so that
the inequalities between this $\lambda$ and the first non-zero eigenvalue
of $\Delta$ are reversed.)

In the next proposition and for the rest of this section, $\Gamma^j$
denotes the $j$-subdivision of (the edges of) the graph $\Gamma$.

\begin{prop}
\label{prop:random SC}
 For any $v\in\naturals$, any $\alpha>0$ and $\alpha'>0$, for any $C\geq 1$, there exists an
integer $j_0$ such that for any $j\geq j_0$, for any graph
$\Gamma$ satisfying the conditions:
\begin{enumerate}
\item Each vertex of $\Gamma$ is of valence at most $v$;
\item The girth of $\Gamma$ is $g$;
\item $\diam(\Gamma)\leq Cg$ for all $i$;
\end{enumerate}
then the following properties hold with probability tending to~$1$
as $g\rightarrow\infty$:
\begin{enumerate}

\item The folded graph $\overline{\Gamma^j}$ obtained by a random
labelling of $\Gamma^j$ satisfies the $\Gr(\alpha)$ condition.

\item There is a reduced word of length between $1$ and
$\alpha'\girth{\overline{\Gamma^j}}$ not appearing on any path in
$\overline{\Gamma^j}$.


\end{enumerate}
\end{prop}

This will be proven in the next sections (a sketch of proof can
also be found in~\cite{Gromov2003}).

Let us now just gather propositions~\ref{prop:expander fact},
\ref{prop:random T} and~\ref{prop:random SC}.

\begin{proof}[Proof of Proposition~\ref{prop:sixth T W}]

Let $\alpha$ be the small cancellation constant to be achieved.

Apply Proposition~\ref{prop:expander fact} with some $v\in \naturals$ to get an
infinite family of graphs $\Gamma_i$; let $\lambda_0$ be the lower
bound on the spectral gap so obtained, and let $C$ be as in this
proposition.  Let us denote by $\Gamma_{i(g)}$ the first graph in this
family having girth at least $g$.

For the chosen $\alpha>0$, let $j$ and $g$ be large enough for the
conclusions of Proposition~\ref{prop:random SC} to hold when
applied to $\Gamma_{i(g)}$. Let $g$ be still large enough
(depending on $j$) so that the conclusions of
Proposition~\ref{prop:random T} applied to this $j$ hold. This
provides a graph satisfying the three requirements of
Proposition~\ref{prop:sixth T W}.
\end{proof}

%
%

\subsection{Some simple properties of random words} Recall $m\geq 2$ is
the number of generators we use. We denote by $\norm{w}$ the norm
in the free group of the word $w$, that is, the length of the associated
reduced word. 

Hereafter $\theta$ is the gross cogrowth of the free group (we refer to
the paragraph ``Growth, cogrowth, and gross cogrowth'' in~\cite{Ollivier04} for
basic properties). Basically, $\theta$ is the infimum of the real numbers so
that the number of words of length $\ell$ which freely reduce to the
trivial word is at most $(2m)^{\theta\ell}$ for all $\ell\in\naturals$.
In particular, the probability that a random walk in the
free group comes back at its origin at time $\ell$ is at most
$(2m)^{-(1-\theta)\ell}$. Explicitly we have $(2m)^{\theta}=2\sqrt{2m-1}$
\cite{Kesten59}.

We state here some elementary properties having to deal with the behavior of
reducing a random word. The first one is pretty intuitive.

\begin{lem}
\label{lem:random reduced}
Let $W_\ell$ be a random word of length $\ell$ and let $\overline{W}_\ell$ be the
associated reduced word. Then the law of $\overline{W}_\ell$ knowing its
length $\abs{\overline{W}_\ell}=\norm{W_\ell}$
is the uniform law on all reduced words of this length.
\end{lem}

\begin{proof}[Proof of the lemma]
The group of automorphisms of the $2m$-regular tree preserving some
basepoint acts transitively on the points at a given distance from the
basepoint and preserves the law of the random walk beginning at this basepoint.
\end{proof}

The following is proven in \cite{Ollivier04}, Proposition~17.

\begin{lem}
\label{lem:random norm}
Let $W_\ell$ be a random word of length $\ell$. Then, for any
$0\leq L\leq \ell$ we have
\[
\Pr(\norm{W_\ell} \leq L ) \leq
(2m)^{-\ell(1-\theta)+\theta L}
\]
\end{lem}

Note that exponent vanishes for $L=\frac{1-\theta}{\theta}\ell<\ell$
(since $\theta>1/2$).
A slightly different, asymptotically stronger version of this lemma
is the following.

\begin{lem}
\label{lem:other random norm}
Let $W_\ell$ be a random word of length $\ell$. Then, for any $L$ we have
\[
\Pr(\norm{W_\ell} \leq L)\leq
\sqrt{\ell\,\frac{2m}{2m-1}}\,(2m)^{-(1-\theta)\ell}(2m-1)^{L/2}
\]
\end{lem}

\begin{proof}
Let $B_\ell$ be the ball of radius $\ell$ centered at $e$ in the free group.
Let $p^\ell_x$ be the probability that $W_\ell=x$.
We have
\begin{eqnarray*}
\mathbb{E}(2m-1)^{-\frac12\norm{W_\ell}}
&=& \sum_{x\in B_\ell} p^\ell_x\, (2m-1)^{-\frac12\norm{x}}
\\
&\leq& \sqrt{\sum_{x\in B_\ell}(p^\ell_x)^2}\,
\sqrt{\sum_{x\in B_\ell} (2m-1)^{-\norm{x}}}
\end{eqnarray*}
by the Cauchy-Schwarz inequality. But $\sum_{x\in B_\ell}(p^\ell_x)^2$
is exactly the probability of return to $e$ at time $2\ell$ of the random
walk (condition by where it is at time $\ell$) which is at most
$(2m)^{-2(1-\theta)\ell}$. Besides, there are $(2m)(2m-1)^{k-1}$ elements
of norm $k$ in $B_\ell$, so that $\sum_{x\in B_\ell} (2m-1)^{-\norm{x}}=
\sum_{0\leq k\leq \ell} (2m)(2m-1)^{k-1}(2m-1)^{-k}=\ell\frac{2m}{2m-1}$.
So we get
\[
\mathbb{E}(2m-1)^{-\frac12\norm{W_\ell}}
\leq \sqrt{\ell\,\frac{2m}{2m-1}} \,(2m)^{-(1-\theta)\ell}
\]

Now we simply apply the Markov inequality
\begin{align*}
\Pr(\norm{W_\ell}\leq L)
&=
\Pr\left((2m-1)^{-\frac12\norm{W_\ell}}\geq (2m-1)^{-\frac12 L}\right)
\\&\leq (2m-1)^{\frac12 L}\,\mathbb{E}(2m-1)^{-\frac12 \norm{W_\ell}}
\end{align*}
to get the conclusion.
\end{proof}

\subsection{Folding the labelled graph}

Labelling a graph by plain random words does generally not result
in a reduced labelling. Nevertheless, we can always fold the
resulting labelled graph. Here we show that in the circumstances
needed for our applications, this folding is a quasi-isometry.
This will allow a transfer of the $\Gr$ small cancellation
condition from the unfolded to the folded graph.

\begin{prop}
For any $\beta>0$, for any $v\in\naturals$, for any $C\geq 1$, there exists an integer
$j_0$ such that for any $j\geq j_0$, for any graph $\Gamma$
satisfying the conditions:
\begin{enumerate}
\item Vertices of $\Gamma$ are of valency at most $v$.
\item $\diam(\Gamma)\leq Cg$ for all $i$, where $g$ is the girth of
$\Gamma$.
\end{enumerate}
then the folding map $\Gamma^j\rightarrow \overline{\Gamma^j}$
from a random labelling $\Gamma^j\rightarrow B$ to the associated
reduced labelling $\overline{\Gamma^j}\looparrowright B$ is a
$(\frac{\theta}{1-\theta},\beta j g,gj)$ local quasi-isometry,
with probability tending to $1$ as $g\rightarrow\infty$.
\end{prop}

We use the notation from~\cite{GhysBook90} for local quasi-isometries:
an $(a,b,c)$ local quasi-isometry is a map $f$ such that whenever
$d(x,y)\leq c$ we have $\frac1a d(x,y)-b\leq d(f(x),f(y))\leq a
d(x,y)+b$. Here folding obviously decreases distances so that only the
left inequality has to be checked.

\begin{rem}
\label{rem:numberofpaths}
Below we will make repeated use of the following: The number of
paths of length $\ell$ in $\Gamma^j$ is at most
$j^2\,v^{Cg+\ell/j}$. Indeed, the number of points in $\Gamma$ is at
most $v^{Cg}$, and once a point is chosen the number of paths of
length $k$ originating at it is at most $v^k$. Now
specifying a path in the subdivision $\Gamma^j$ amounts to
specifying a path in $\Gamma$ and giving two integers between $1$
and $j$ to specify the exact endpoints.
\end{rem}

\begin{proof}
Unwinding the definition of local quasi-isometries, we have to prove that any immersed path of length $\beta g j +
\ell\leq gj$ in $\Gamma^j$ is mapped onto a path of length at
least $\frac{1-\theta}{\theta}\ell$ in $\overline{\Gamma^j}$.

By Remark~\ref{rem:numberofpaths}, there are at most $j^2\,v^{Cg+g}$ paths of length
$gj$ in the subdivision $\Gamma^j$ of $\Gamma$. Fix such a path, of
length say $\beta g j + \ell$.

Since the length of the immersed path is at most $gj=\girth(\Gamma^j)$, the path
does not travel twice along the same edge. Consequently, the
labels appearing on this path are all chosen independently. Then
by Lemma~\ref{lem:random norm}, the probability that its length after
folding is less than $\frac{1-\theta}{\theta}\ell$ is less than
\[
(2m)^{-(1-\theta)(\ell+\beta g j)+\theta
\frac{1-\theta}{\theta}\ell}=(2m)^{-(1-\theta)\beta g j}
\]
for this particular path. Since the number of choices for the path
is at most $j^2\,v^{Cg+g}$, if $j$ is large enough depending on
$C$, $\beta$ and $\theta$, namely if
$v^{C+1}(2m)^{-(1-\theta)\beta j} < 1$, then the probability that
\emph{there exists} a path violating our local quasi-isometry
property will tend to $0$ as $g\rightarrow \infty$.
\end{proof}

\begin{cor}
\label{cor:folded girth}
In the same circumstances, the girth of $\overline{\Gamma^j}$ is
at least $\frac{1-\theta}{\theta}-\beta$ times that of $\Gamma^j$.
\end{cor}

\begin{proof}
Take a simple closed path $p$ in $\overline{\Gamma_j}$. It is the
image of a non-null-homotopic closed path $q$ in $\Gamma^j$, whose
length is by definition at least $gj=\girth \Gamma^j$. Let $q'$ be
the initial subpath of $q$ of length $gj$. We can apply the local
quasi-isometry statement to $q'$, showing that its image $p'$ has
length at least $\frac{1-\theta}{\theta}gj-\beta gj$, which is
thus a lower bound on the length of $p$.
\end{proof}

\subsection{Pieces in the unfolded and folded graphs.} Here we
show that under the circumstances above, the probability to get a
long piece in the folded graph is very small.

Suppose again that we are given a graph $\Gamma$ of degree at most
$v$, of girth $g$ and of diameter at most $Cg$. Consider its
$j$-subdivision $\Gamma^j$ endowed with a random labelling and let
$\overline{\Gamma^j}$ be the associated folded labelled graph.

Let $p, p'$ be two immersed paths in $\overline{\Gamma^j}$. Let $q,q'$ be
some preimages in $\Gamma^i$ of $p,p'$. If $p$ and $p'$ are labelled by
the same word, then $q$ and $q'$ will be labelled by some freely equal
words, so that pieces come from pieces.

Note that in a graph labelled by non-reduced words, there are some
``trivial pieces'': e.g.\ if some $aa^{-1}$ appears next to a word $w$,
then $(w,aa^{-1}w)$ will be a piece. Such pieces disappear after folding
the labelled graph; this is why we discard them in the following.

\begin{prop}
\label{prop:random pieces}
Let $q$, $q'$ be two immersed paths in a graph $\Delta$ of girth
$g$. Suppose that $q$ and $q'$ have length $\ell$ and $\ell'$
respectively, with $\ell$ and $\ell'$ at most $g/2$.
Endow $\Delta$ with a random labelling. Suppose that after folding
the graph, the paths $q$ and $q'$ are mapped to distinct paths.
Then the probability that $q$ and $q'$ are labelled by two freely
equal words is at most
\[
C_{\ell,\ell'}(2m)^{-(1-\theta)(\ell+\ell')}
\]
where $C_{\ell,\ell'}$ is a term growing subexponentially in $\ell+\ell'$.
\end{prop}

\begin{proof}
Let $w$ and $w'$ be the words labelling $q$ and $q'$ respectively.

First, assume that the images of $q$ and $q'$ in $\Delta$ are
disjoint. Then the letters making up $w$ and $w'$ are chosen
independently, and thus the word $w{w'}^{-1}$ is a plain random
word. Thus is this case the proposition is just a rewriting of the
definition of $\theta$.

Second, suppose that the paths do intersect in $\Delta$: this results in
lack of independence in the choice of the letters making up $w$ and $w'$
(the same problem is treated in a slightly different setting
in~\cite{Ollivier04}, section ``Elimination of doublets''), which needs
to be treated carefully.  Since the length of these words is less than
half the girth, the intersection in $\Delta$ is connected and we can
write $w=u_1u_2u_3$, $w'=u'_1u_2u'_3$ where the $u_i$'s are
\emph{independently chosen} random words (depending on relative
orientation of $w$ and $w'$, $u_2^{-1}$ rather than $u_2$ may appear in
$w'$). We can suppose that $u'_1u_1^{-1}$ is not freely trivial:
otherwise the two paths start at the same point after folding, and so if
$w=w'$ we also have $u'_3u_3^{-1}=e$ so that they also end at the same
point after folding, but this is discarded by assumption.  Likewise
$u'_3u_3^{-1}$ is not freely trivial.

Let $v_1$, $v_2$, $v_3$, $v'_1$, $v'_3$ be the reduced words freely equal
to $u_1$, $u_2$\ldots respectively.

Lemma~\ref{lem:random reduced} tells us that the words $v_1$, $v_2$\ldots are random reduced words.
Now let us draw a picture expressing the equality $v_1v_2v_3=v'_1v_2v'_3$:
\begin{figure}\centering
\includegraphics[width=.5\textwidth]{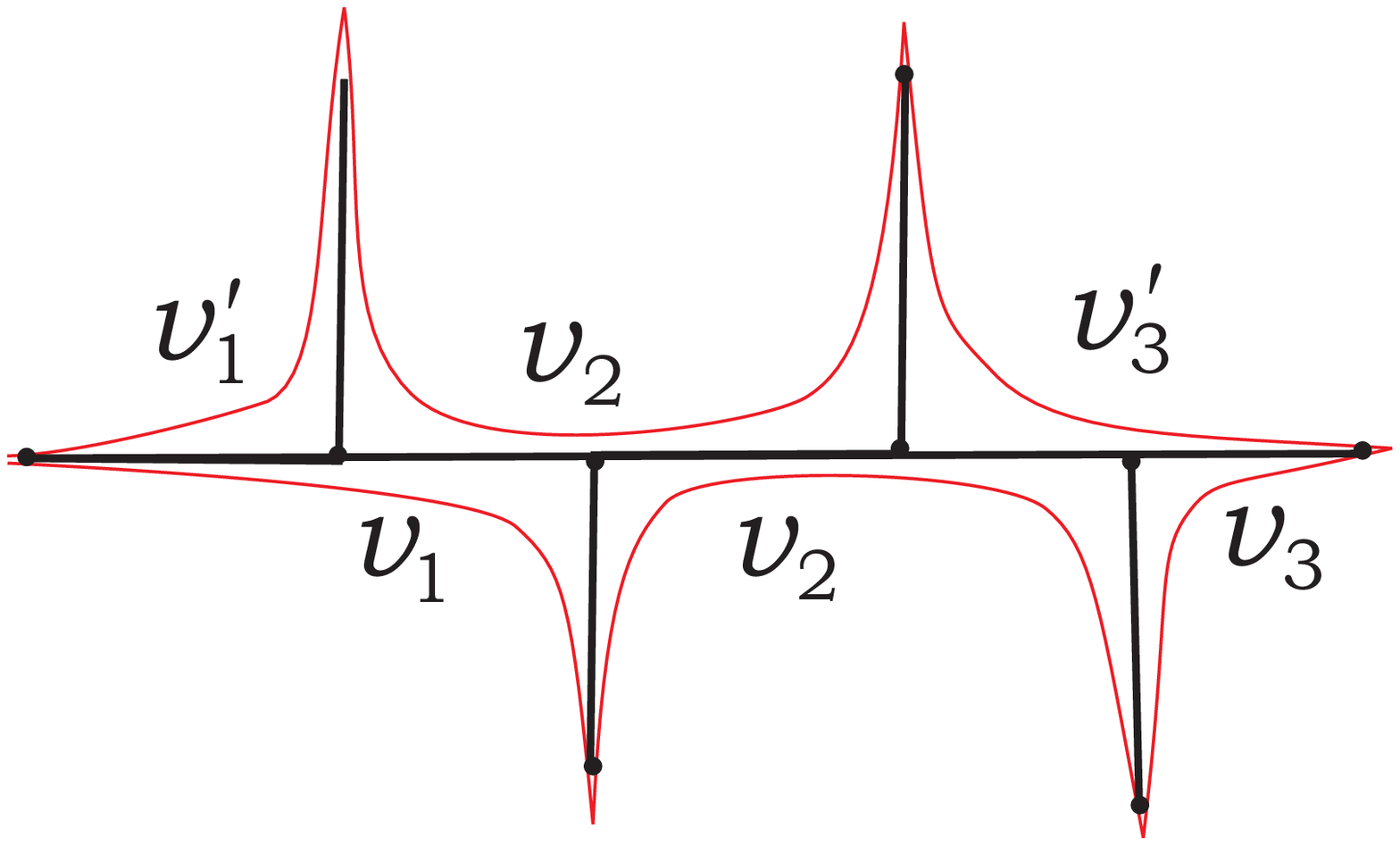}
\caption{ \label{fig:foldedhexagon}}
\end{figure}

Note that the two copies of $v_2$ have to be shifted relatively to each other,
otherwise this means that $u'_1u_1^{-1}$ and $u'_3u_3^{-1}$ are freely
trivial.

Let $k$ be the length shared between the two copies of $v_2$.  Now let us
evaluate the probability of this situation knowing all the lengths of the
words $v_1$, $v_2$, \ldots Conditionnally to their lengths, these words
are uniformly chosen random reduced words by Lemma~\ref{lem:random
reduced}.

We begin with the two copies of $v_2$: though they are not chosen
independently, since we know that they are shifted, adding letter after
letter we see that the probability that they can glue along a subpath of
length $k$ is a most $1/(2m-1)^k$. Once $v_2$ is given, the words $v_1$,
$v_3$, $v'_1$, $v'_3$ are all chosen independently of each other. The
probability that they glue according to the picture is $1/(2m-1)^{L-k}$
where $L$ is the total length of the picture. So the overall probability
of such a gluing is $1/(2m-1)^L$.

We obviously have $\ell+\ell'=|v_1|+2|v_2|+|v_3|+|v'_1|+|v'_3|=2L$.  Now
by Proposition~\ref{lem:other random norm} applied to all these words
separately, the probability of achieving this value of
$|v_1|+2|v_2|+|v_3|+|v'_1|+|v'_3|$ is less than
\[
C_{\ell,\ell'}\,(2m)^{-(1-\theta)(|u_1|+2|u_2|+|u_3|+|u'_1|+|u'_3|)}(2m-1)^{\frac12\,2L}=
C_{\ell,\ell'}\,(2m)^{-(1-\theta)(\ell+\ell')}(2m-1)^L
\]
where $C_{\ell,\ell'}$ is a term growing subexponentially in $\ell+\ell'$.

So the overall probability of such a situation, taking into account the
possibilities for $L$ between $0$ and $\ell+\ell'$, is at most
\[
\sum_{0\leq L \leq
\ell+\ell'}(2m-1)^{-L}\,C_{\ell+\ell'}\,(2m)^{-(1-\theta)(\ell+\ell')}(2m-1)^L
\]
since we just proved above that $(2m-1)^{-L}$ is an upper bound for the probability
of the situation knowing $L$.  
But this is equal to $C'_{\ell+\ell'}\,(2m)^{-(1-\theta)(\ell+\ell')}$
where $C'_{\ell+\ell'}$ is another term growing subexponentially in
$\ell+\ell'$.
\end{proof}

We are now ready to prove Proposition~\ref{prop:random SC} stating that
the $\Gr(\alpha)$ condition is satisfied with overwhelming probability.
In order to avoid heavy notations, we will directly prove the stronger
variant of the $\Gr$ condition involving the girth instead of the length
of cycles containing the pieces (see section~\ref{section:graphical}).

\begin{proof}[Proof of Proposition~\ref{prop:random SC}, small cancellation part]
Since ruling out small pieces rules out larger pieces as well, it is
enough to work for small $\alpha$.

Let $\bar{g}$ be the girth of $\overline{\Gamma^j}$.  By
Corollary~\ref{cor:folded girth}, we can assume that $\bar{g}\geq
\left(\frac{1-\theta}{\theta}-\beta\right)gj$ with overwhelming
probability, for arbitrarily small $\beta$.

Let $p,p'$ be two distinct immersed paths in $\overline{\Gamma^j}$
forming a $\alpha$-piece; both $p$ and $p'$ are of length $\alpha
\bar{g}$.  Let $q$ and $q'$ be some immersed paths in $\Gamma^j$ mapping
to $p$ and $p'$.

Suppose that the length of $q$ (or $q'$) is greater than $gj/2$. By
applying the local quasi-isometry property to an initial subpath of $q$
of length $gj/2$ we get that the length of $p$ would be at least
$\frac{1-\theta}{\theta}gj/2-\beta gj$. But the length of $p$ is exactly
$\alpha\bar{g}\leq \alpha gj$, so that if $\alpha$ and $\beta$ are taken
small enough (depending on $\theta$) we get a contradiction. Hence, the
length of $q$ is at most $gj/2$, so that we are in a position to apply
Proposition~\ref{prop:random pieces}.

The length of $q$ and $q'$ is at least that of $p$ and $p'$ namely
$\alpha \bar{g}$, and since $\bar{g}\geq
\left(\frac{1-\theta}{\theta}-\beta\right)gj$, $q$ and $q'$ form a
$\alpha\left(\frac{1-\theta}{\theta}-\beta\right)$-piece in $\Gamma^j$.
Now Proposition~\ref{prop:random pieces} states that for fixed $q$
and $q'$ in $\Gamma^j$, the probability of this is at most
$C_{gj}(2m)^{-(1-\theta)2gj\alpha\left(\frac{1-\theta}{\theta}-\beta\right)}$,
where $C_{gj}$ is a subexponential term in $\abs{q}+\abs{q'}\leq gj$.

By Remark~\ref{rem:numberofpaths}, the number of choices for $q$ and $q'$ is at
most $j^4\,v^{(2C+1)g}$.
So the probability that one of these choices gives rise to a
piece is at most
\[
j^4\,v^{(2C+1)g}\,C_{gj}(2m)^{-(1-\theta)2gj\alpha\left(\frac{1-\theta}{\theta}-\beta\right)}
\]

Now, if $\beta$ is taken small enough (depending only on $\theta$)
and if $j$ is taken large enough (depending on $\alpha, \theta$
and $C$ but not on $g$), namely if
\[
v^{2C+1}\,(2m)^{-(1-\theta)2j\alpha\left(\frac{1-\theta}{\theta}-\beta\right)}<1
\]
then this tends to $0$ when $g$ tends to infinity.
\end{proof}

\begin{proof}[Proof of Proposition~\ref{prop:random SC}, missing
word part] We now prove that for any $\alpha'>0$, in the same
circumstances, there exists a reduced word of length $\alpha'
\girth(\overline{\Gamma^j})$ not appearing on any path in
$\overline{\Gamma^j}$.

Let $p$ be a simple path of length $\alpha' \bar{g}$ in
$\overline{\Gamma^j}$. It is the image of some path
$q$ in $\Gamma^j$ of length at least $\alpha' \bar{g}\geq
\alpha'(\frac{1-\theta}{\theta}-\beta)gj$. But by
Remark~\ref{rem:numberofpaths}, the
number of such paths in $\Gamma^j$ is at most
$j^2\,v^{Cg+\alpha'(\frac{1-\theta}{\theta}-\beta)g}$, whereas the
total number of reduced words of this length is at least
$(2m-1)^{\alpha'(\frac{1-\theta}{\theta}-\beta)gj}$. So if $j$ is
taken large enough (depending on $\alpha'$ and $\theta$ but not on $g$) that is if
\[
v^{C+\alpha'(\frac{1-\theta}{\theta}-\beta)}<(2m-1)^{\alpha'(\frac{1-\theta}{\theta}-\beta)j}
\]
then the possible reduced words outnumber the paths in
$\overline{\Gamma^j}$ when $g\rightarrow\infty$, so that there has to be
a missing word.
\end{proof}

\section{Problems} \begin{com}
Probably better to integrate into paper, but might as well collect
them here for now... \end{com}
 Does there exist a finitely
presented group $N$ with property~$T$ such that $\Out(N)$ is
infinite?

Let $Q$ be a f.p.\ group with property~$T$. Does there exist
word-hyperbolic $G$ with property~$T$
 and f.g.\ normal subgroup $N$ such that
 $Q=G/N$?

Do there exist f.p.\ Kazhdan groups which are not Hopfian or coHopfian?

\bibliographystyle{alpha}

\bibliography{wise}

%
%
\end{document}